\documentclass[11pt,a4paper]{article}

\usepackage{graphicx}
\usepackage{latexsym}
\usepackage{eepic}
\usepackage{amsfonts}

\makeatletter
\newdimen\mynormalparindent
\def\mymakefnmark{}
\def\mymakefntext{\indent\mymakefnmark}
\long\def\myfootnotetext#1{\insert\footins{%
  \normalfont\footnotesize
  \interlinepenalty\interfootnotelinepenalty
  \splittopskip\footnotesep \splitmaxdepth \dp\strutbox
  \floatingpenalty\@MM \hsize\columnwidth
  \@parboxrestore \parindent\mynormalparindent \sloppy
  \mymakefntext{\rule\z@\footnotesep\ignorespaces#1\unskip\strut\par}}}

\makeatother

\long\def\symbolfootnote[#1]#2{\begingroup%
\def\thefootnote{\fnsymbol{footnote}}\footnote[#1]{#2}\endgroup}

\newtheorem{basic}{Basic}[section]

\newtheorem{lem}[basic]{Lemma}
\newtheorem{propos}[basic]{Proposition}
\newtheorem{thm}[basic]{Theorem}

\newcommand{\bdm}{\begin{displaymath}}
\newcommand{\edm}{\end{displaymath}}
\newcommand{\be}{\begin{equation}}
\newcommand{\ee}{\end{equation}}
\newcommand{\ep}{\vspace{-3mm}\hfill\mbox{$\Box$}\\}
\newcommand{\ul}{\underline}

\newcommand{\R}{\mathbb{R}}

\begin{document}

\begin{center}
{\bf\Large The Generalized Schur Decomposition and the\\ rank-$R$ set of
real $I\times J\times 2$ arrays}\\

\vspace{2cm} Alwin Stegeman
\symbolfootnote[2]{The author is with the Heijmans Institute for
Psychological Research, University of Groningen, Grote Kruisstraat 2/1,
9712 TS Groningen, The Netherlands, phone: ++31 50 363 6193, fax: ++31 50
363 6304, email: a.w.stegeman@rug.nl, URL:
http://www.gmw.rug.nl/$\sim$stegeman. The author is supported by the
Dutch Organisation for Scientific Research (NWO), VIDI grant 452-08-001.}\\
\vspace{5mm} \today \vspace{1cm}

\begin{abstract}
\noindent It is known that a best low-rank approximation to multi-way
arrays or higher-order tensors may not exist. This is due to the fact that
the set of multi-way arrays with rank at most $R$ is not closed.
Nonexistence of the best low-rank approximation results in diverging
rank-1 components when an attempt is made to compute the approximation.
Recently, a solution to this problem has been proposed for real $I\times
J\times 2$ arrays. Instead of a best rank-$R$ approximation the best
fitting Generalized Schur Decomposition (GSD) is computed. Under the
restriction of nonsingular upper triangular matrices in the GSD,
the set of GSD solutions equals the interior and boundary of the rank-$R$
set. Here, we show that this holds even without the restriction.
We provide a complete classification of interior, boundary, and exterior
points of the rank-$R$ set of real $I\times J\times 2$ arrays, and show
that the set of GSD solutions equals the interior and boundary of this set.

\end{abstract}\end{center}

\newpage
\section{Introduction}
\setcounter{equation}{0}
This paper is an addendum to Stegeman and De Lathauwer \cite{SDL} who study
the following subject. Let $\circ$ denote the outer-product, and define
the outer-product rank of $\ul{\bf Y}\in \R^{I\times J\times K}$ as
\be
{\rm rank}_{\circ}(\ul{\bf Y})=\min\{R\;|\;\ul{\bf Y}=\sum_{r=1}^R
{\bf x}_r\circ{\bf y}_r\circ{\bf z}_r\}\,.
\ee

\noindent Let
\be
{\cal S}_R(I,J,K)=\{\ul{\bf Y}\in\R^{I\times J\times K}\;|\;
{\rm rank}_{\circ}(\ul{\bf Y})\le R\}\,,
\ee

\noindent and let $\overline{\cal S}_R(I,J,K)$ denote the closure of ${\cal
S}_R(I,J,K)$, i.e. the union of the set itself and its boundary points in
$\R^{I\times J\times K}$.

Let $\ul{\bf Z}\in\R^{I\times J\times K}$ and $||\cdot||$ denote the
Frobenius norm on $\R^{I\times J\times K}$. Consider the following low-rank
approximation problem. \be \label{prob-CP} \min\{||\ul{\bf Z}-\ul{\bf
Y}||\;|\;\ul{\bf Y}\in{\cal S}_R(I,J,K)\}\,. \ee

\noindent Assuming rank$_{\circ}(\ul{\bf Z})>R$, an optimal solution
of (\ref{prob-CP}) will be a boundary point of the set ${\cal
S}_R(I,J,K)$. However, the set $S_R(I,J,K)$ is not closed for $R\ge 2$, and
problem (\ref{prob-CP}) may not have an optimal solution due to this fact;
see De Silva and Lim \cite{DSL}. Nonexistence of an optimal solution results in
diverging rank-1 components when an attempt is made to compute a best
rank-$R$ approximation, see Krijnen, Dijkstra and Stegeman \cite{KDS}.
In order to overcome this fallacy, \cite{DSL} proposed to consider
instead \be \label{prob-closure} \min\{||\ul{\bf Z}-\ul{\bf Y}||\;|\;\ul{\bf
Y}\in\overline{\cal S}_R(I,J,K)\}\,. \ee

\noindent Note that if (\ref{prob-CP}) has an optimal solution, then
it is also an optimal solution of (\ref{prob-closure}). To solve
problem (\ref{prob-closure}), we need to characterize the boundary points
of ${\cal S}_R(I,J,K)$ and we need an algorithm to find an optimal boundary
point. For $R=2$, the boundary points are determined in \cite{DSL}, and an
algorithm to solve (\ref{prob-closure}) is proposed in Rocci and Giordani
\cite{RoG}. For a general approach to obtain an optimal solution to
(\ref{prob-closure}) from an attempt to solve (\ref{prob-CP}), see Stegeman
\cite{Ste-divcomp}.

In Stegeman and De Lathauwer \cite{SDL} the case $K=2$ is considered. Let
\newpage
\bdm
{\cal P}_R(I,J,2)=\{\ul{\bf Y}\in\R^{I\times J\times 2}\;|\;
{\bf Y}_k={\bf Q}_a\,{\bf R}_k\,{\bf Q}_b^T,\,k=1,2,\;{\rm with}\;
{\bf Q}_a^T{\bf Q}_a={\bf Q}_b^T{\bf Q}_b={\bf I}_R
\edm
\be
\label{eq-GSDset}
\hspace{2cm}{\rm and}\;{\bf R}_k\;(R\times R)\;{\rm upper\;triangular}\}\,,
\ee

\noindent denote the set of arrays with a full Generalized Schur
Decomposition (GSD). Here, ${\bf Y}_k$ ($I\times J$) denotes the $k$th
frontal slice of $\ul{\bf Y}$. Note that a GSD exists only for
$R\le\min(I,J)$. In \cite{SDL} it is shown that the problem \be
\label{prob-GSD} \min\{||\ul{\bf Z}-\ul{\bf Y}||\;|\;\ul{\bf Y}\in{\cal
P}_R(I,J,2)\}\,, \ee

\noindent is guaranteed to have an optimal solution. Moreover, it holds that
${\cal P}_R(I,J,2)=\overline{\cal S}_R(I,J,2)$ under the restriction that
only arrays are considered that have a GSD with ${\bf R}_1$ and ${\bf R}_2$
nonsingular. Also, a Jacobi algorithm (based on De Lathauwer, De Moor and
Vandewalle \cite{LMV}) is presented for solving (\ref{prob-GSD}). Hence,
under the above restriction, for $K=2$ problem (\ref{prob-closure}) can
be solved by solving problem (\ref{prob-GSD}).

In this note we show that the restriction used in
\cite{SDL} is not necessary. That is, we prove that ${\cal
P}_R(I,J,2)=\overline{\cal S}_R(I,J,2)$ holds for $R\le\min(I,J)$.

We use the notation $({\bf S},{\bf T},{\bf U})\cdot\ul{\bf Y}$ to denote
the multilinear matrix multiplication of an array $\ul{\bf Y}\in\R^{I\times
J\times K}$ with matrices ${\bf S}$ ($I_2\times I$), ${\bf T}$ ($J_2\times
J$), and ${\bf U}$ ($K_2\times K$). The result of the multiplication is an
$I_2\times J_2\times K_2$ array. We refer to $({\bf I}_I,{\bf I}_J,{\bf
U})\cdot\ul{\bf Y}$ with ${\bf U}$ ($K\times K$) nonsingular as a slicemix.

For later use we mention that, for nonsingular ${\bf S}$, ${\bf T}$, ${\bf
U}$, and $\ul{\bf X}=({\bf S},{\bf T},{\bf U})\cdot\ul{\bf Y}$, we have
rank$_{\circ}(\ul{\bf X})=$ rank$_{\circ}(\ul{\bf Y})$ and $\ul{\bf X}$
is an interior (boundary, exterior) point of ${\cal S}_R(I,J,K)$ if and
only if $\ul{\bf Y}$ is an interior (boundary, exterior) point of ${\cal
S}_R(I,J,K)$.

\section{The case $I=J=R$}
\setcounter{equation}{0}
Here, we consider the case where the arrays have two $I\times I$ slices and
the number of components equals $I$. In \cite{SDL} only arrays are
considered that have two nonsingular slices. In Proposition~\ref{p-1} below
we present a complete classification of $I\times I\times 2$ arrays into
interior, boundary, and exterior points of the set ${\cal S}_R(I,J,2)$.
This classification is used to show that ${\cal P}_I(I,I,2)=\overline{\cal
S}_I(I,I,2)$ in Theorem~\ref{t-1}. In the proofs of Proposition~\ref{p-1}
and Theorem~\ref{t-1} we use the following lemma which extends the
theory on real matrix pencils and may be of interest in itself.

\begin{lem}
\label{lem-QZ}
Let $\ul{\bf Y}\in\R^{I\times I\times 2}$ with $I\times I$ slices
${\bf Y}_1$ and ${\bf Y}_2$ such that {\rm det}$(\mu\,{\bf
Y}_1+\lambda\,{\bf Y}_2)=0$ for all $\mu,\lambda\in\R$. Then
$\ul{\bf Y}\in{\cal P}_I(I,I,2)$.
\end{lem}

\noindent {\bf Proof.} As shown in Moler and Stewart \cite{MoS} (see also
Golub and Van Loan \cite[Section 7.7.2]{matcomp}) there exist orthonormal
${\bf Q}$ and ${\bf Z}$ such that ${\bf G}={\bf Q}\,{\bf Y}_2\,{\bf Z}$ is
upper triangular and ${\bf F}={\bf Q}\,{\bf Y}_1\,{\bf Z}$ is quasi-upper
triangular. That is, ${\bf F}$ is block-upper triangular where its diagonal
blocks are $2\times 2$ or $1\times 1$ in size. The proof is complete if we
show that there exist orthonormal $\tilde{\bf Q}$ and $\tilde{\bf Z}$ such
that $\tilde{\bf Q}\,{\bf F}\,\tilde{\bf Z}$ and $\tilde{\bf Q}\,{\bf
G}\,\tilde{\bf Z}$ are upper triangular.

Let ${\bf F}$ have diagonal blocks ${\bf F}_1,\ldots,{\bf F}_m$, where ${\bf
F}_i$ is $2\times 2$ or $1\times 1$. Denote the corresponding diagonal
blocks of ${\bf G}$ by ${\bf G}_i$, $i=1,\ldots,m$. We have
\be
{\rm det}(\mu\,{\bf Y}_1+\lambda\,{\bf Y}_2)=\prod_{i=1}^m
{\rm det}(\mu\,{\bf F}_i+\lambda\,{\bf G}_i)=0\,,\quad\quad\quad
{\rm for\;all\;}\mu,\lambda\in\R\,.
\ee

\noindent This can only hold if
\be
\label{eq-oneblock}
{\rm for\;some\;} l\,,\quad\quad {\rm det}(\mu\,{\bf F}_l+\lambda\,{\bf
G}_l)=0\,,\quad\quad\quad {\rm for\;all\;} \mu,\lambda\in\R\,.
\ee

\noindent It is shown in Moler and Stewart \cite[Section 5]{MoS} that if
${\bf F}_i$ is a $2\times 2$ block and det$({\bf F}_i+\lambda\,{\bf
G}_i)=0$ for some $\lambda\in\R$, then $2\times 2$ orthonormal $\tilde{\bf
Q}$ and $\tilde{\bf Z}$ can be found such that $\tilde{\bf Q}\,{\bf
F}_i\,\tilde{\bf Z}$ and $\tilde{\bf Q}\,{\bf G}_i\,\tilde{\bf Z}$ are
upper triangular. Hence, we may assume without loss of generality that if
${\bf F}_i$ is a $2\times 2$ block, then det$({\bf F}_i+\lambda\,{\bf
G}_i)\neq 0$ for all $\lambda\in\R$.

Let index $l$ be as in (\ref{eq-oneblock}). From the discussion above it
follows that we may assume that ${\bf F}_l$ is $1\times 1$. Hence, ${\bf
F}$ and ${\bf G}$ have a zero on their diagonals in the same position.
Suppose the common zero appears right after a $2\times 2$ block ${\bf F}_i$,
i.e.
\be
\label{eq-33mat}
\left[\begin{array}{cc}
{\bf F}_i & {\bf f} \\
{\bf 0}^T & 0 \end{array}\right]
=\left[\begin{array}{cc|c}
* & * & * \\
~* & * & * \\
\hline
0 & 0 & 0 \end{array}\right]\,,\quad\quad\quad\quad\quad
\left[\begin{array}{cc}
{\bf G}_i & {\bf g} \\
{\bf 0}^T & 0\end{array}\right]
=\left[\begin{array}{cc|c}
* & * & * \\
0 & * & * \\
\hline
0 & 0 & 0 \end{array}\right]\,.
\ee

\noindent Let ${\bf z}_1\in\R^3$ be orthogonal to the second rows of the
matrices in (\ref{eq-33mat}). Then postmultiplying (\ref{eq-33mat}) by any
orthonormal $\tilde{\bf Z}=[{\bf z}_1|{\bf z}_2|{\bf z}_3]$ brings both
matrices into upper triangular form and leaves the common zero in the same
position. Analogously, if the common zero appears right before a $2\times
2$ block ${\bf F}_i$, then premultiplying by a suitable orthonormal
$\tilde{\bf Q}$ does the trick.

If the common zero on the diagonals of ${\bf F}$ and ${\bf G}$ is not
adjacent to a $2\times 2$ block ${\bf F}_i$, then we resort to
simultaneously reordering the diagonal blocks of ${\bf F}$ and ${\bf
G}$ (except the common zero) such that it is. It suffices to show that
swapping adjacent $2\times 2$ and $1\times 1$ blocks is possible by
orthonormal transformations. Let ${\bf F}_i$ be $2\times 2$ and consider
the $3\times 3$ matrices
\be
\left[\begin{array}{cc}
{\bf F}_i & {\bf f} \\
{\bf 0}^T & f_{i+1} \end{array}\right]
\,,\quad\quad\quad\quad\quad
\left[\begin{array}{cc}
{\bf G}_i & {\bf g} \\
{\bf 0}^T & g_{i+1}\end{array}\right]\,.
\ee

\noindent Swapping the diagonal blocks $i$ and $i+1$ by orthonormal
transformations is possible if ${\bf x},{\bf y}\in\R^2$ exist such that
they satisfy the so-called generalized Sylvester equation (see e.g.
Kressner \cite[Section 5]{Kre}):
\be
\label{eq-gensyl}
{\bf F}_i\,{\bf x} - f_{i+1}\,{\bf y}={\bf f}\,,\quad\quad\quad\quad
{\bf G}_i\,{\bf x} - g_{i+1}\,{\bf y}={\bf g}\,.
\ee

\noindent We may assume that $f_{i+1}$ and $g_{i+1}$ are not both zero. Let
$g_{i+1}\neq 0$ (the proof for $f_{i+1}\neq 0$ is analogous). Then
(\ref{eq-gensyl}) is satisfied for
\be
\label{eq-solxy}
{\bf y}=({\bf G}_i\,{\bf x}-{\bf g})/g_{i+1}\,,\quad\quad\quad\quad
{\bf x}=({\bf F}_i-(f_{i+1}/g_{i+1})\,{\bf G}_i)^{-1}\,
({\bf f}-(f_{i+1}/g_{i+1})\,{\bf g})\,.
\ee

\noindent Note that det$({\bf F}_i-(f_{i+1}/g_{i+1})\,{\bf G}_i)\neq 0$ by
assumption, and the solution (\ref{eq-solxy}) is unique. Hence, the
diagonal blocks can be swapped. Analogously, it can be shown that a
$1\times 1$ block $i$ and a $2\times 2$ block $i+1$ can be swapped. This
completes the proof.\ep

\begin{propos}
\label{p-1}
Let $\ul{\bf Y}\in\R^{I\times I\times 2}$.
\begin{itemize}
\item[{\rm (a)}] If there exists a ${\bf U}$  nonsingular such that
$\ul{\bf X}=({\bf I}_I,{\bf I}_I,{\bf U})\cdot\ul{\bf Y}$ has nonsingular
slice ${\bf X}_1$, then
\begin{itemize}
\item[{\rm (a1)}] $\ul{\bf Y}$ is an interior point of ${\cal S}_I(I,I,2)$
if ${\bf X}_2{\bf X}_1^{-1}$ has $I$ distinct real eigenvalues.
\item[{\rm (a2)}] $\ul{\bf Y}$ is a boundary point of ${\cal S}_I(I,I,2)$
if ${\bf X}_2{\bf X}_1^{-1}$ has $I$ real eigenvalues but not all distinct.
\item[{\rm (a3)}] $\ul{\bf Y}$ is an exterior point of ${\cal S}_I(I,I,2)$
if ${\bf X}_2{\bf X}_1^{-1}$ has at least one pair of complex eigenvalues.
\end{itemize}

\item[{\rm (b)}] If there does not exist a ${\bf U}$ nonsingular
such that $\ul{\bf X}=({\bf I}_I,{\bf I}_I,{\bf U})\cdot\ul{\bf Y}$ has
nonsingular slice ${\bf X}_1$, then $\ul{\bf Y}$ is a boundary
point of ${\cal S}_I(I,I,2)$. \end{itemize}

\end{propos}

\noindent {\bf Proof.} The proofs of (a) follow from the fact that
multilinear matrix multiplication leaves the property interior (boundary,
exterior) point invariant, and application of Stegeman and De
Lathauwer \cite[Lemma 3.1]{SDL} which is due to Stegeman \cite{Ste}.

Next we prove (b). From Lemma~\ref{lem-QZ} it follows that $\ul{\bf
Y}=({\bf Q}_a,{\bf Q}_b,{\bf I}_2)\cdot\ul{\bf R}$, where $\ul{\bf R}$ has
two upper diagonal slices. Below, we show that $\ul{\bf R}$ is a boundary
point of ${\cal S}_I(I,I,2)$. Since ${\bf Q}_a$ and ${\bf Q}_b$
are nonsingular, it follows that also $\ul{\bf Y}$ is a
boundary point of ${\cal S}_I(I,I,2)$.

It holds that det$(\mu\,{\bf R}_1+\lambda\,{\bf R}_2)=0$ for all
$\mu,\lambda\in\R$, which implies that ${\bf R}_1$ and ${\bf R}_2$ have a
zero on their diagonals in the same position. A small perturbation of the
diagonals of ${\bf R}_1$ and ${\bf R}_2$ yields slices
${\bf H}_1$ (nonsingular) and ${\bf H}_2$, with ${\bf H}_2{\bf H}_1^{-1}$
(upper triangular) having $I$ real eigenvalues, and $||\ul{\bf
R}-\ul{\bf H}||<\epsilon$ for any $\epsilon >0$. Next, we show that it is
possible to choose the perturbation such that ${\bf H}_2{\bf H}_1^{-1}$ has
a pair of identical eigenvalues. For simplicity, we assume that the
diagonals of ${\bf R}_1$ and ${\bf R}_2$ contain one common zero. A proof
for the general case is analogous.

Let $({\bf R}_1)_{ii}=({\bf R}_2)_{ii}=0$ and set $({\bf
H}_1)_{ii}=\delta_1$ and $({\bf H}_2)_{ii}=\delta_2$. This yields a nonzero
eigenvalue $\delta_2/\delta_1$ for ${\bf H}_2{\bf H}_1^{-1}$ (assuming
small perturbations of the other zeros on the diagonal of ${\bf R}_1$, such
that ${\bf H}_1$ is nonsingular). Unless stated otherwise, we only perturb
the zero diagonal elements of ${\bf R}_1$ and ${\bf R}_2$. If, for some
$j\neq i$, $({\bf R}_1)_{jj}\neq 0$ and $({\bf R}_2)_{jj}\neq 0$, then let
$\lambda=({\bf R}_2)_{jj}/({\bf R}_1)_{jj}$, and choose
$\delta_2=\lambda\,\delta_1$. This yields ${\bf H}_2{\bf H}_1^{-1}$ with
two identical real eigenvalues $\delta_2/\delta_1=\lambda$ for any
$\delta_1>0$. If no common nonzero diagonal elements of ${\bf R}_1$ and
${\bf R}_2$ exist, then we proceed as follows. If, for some $j\neq i$,
$({\bf R}_1)_{jj}\neq 0$ and $({\bf R}_2)_{jj}=0$, then let $({\bf
H}_2)_{jj}=\eta$, and choose $\delta_1=\sqrt{\delta_2}\,({\bf R}_1)_{jj}$
and $\eta=\sqrt{\delta_2}$. This yields ${\bf H}_2{\bf H}_1^{-1}$ with two
identical real eigenvalues $\delta_2/\delta_1=\eta/({\bf R}_1)_{jj}$ for
any $\delta_2>0$. If, for some $j\neq i$, $({\bf R}_1)_{jj}=0$ and
$({\bf R}_2)_{jj}\neq 0$, then let $({\bf H}_1)_{jj}=\eta$, and choose
$\delta_2=\sqrt{\delta_1}\,({\bf R}_2)_{jj}$ and $\eta=\sqrt{\delta_1}$.
This yields ${\bf H}_2{\bf H}_1^{-1}$ with two identical real eigenvalues
$\delta_2/\delta_1=({\bf R}_2)_{jj}/\eta$ for any $\delta_1>0$. Hence, it
is possible to get ${\bf H}_2{\bf H}_1^{-1}$ with a pair of identical
eigenvalues.

By Proposition~\ref{p-1} (a2), the array $\ul{\bf H}$ is a boundary point
of ${\cal S}_I(I,I,2)$. Since $||\ul{\bf R}-\ul{\bf H}||<\epsilon$ for
any $\epsilon >0$, it follows that $\ul{\bf R}$ can be approximated
arbitrarily closely from $\overline{\cal S}_I(I,I,2)$. Hence, we obtain
$\ul{\bf R}\in\overline{\cal S}_I(I,I,2)$. Moreover, since for any
$\epsilon >0$ the array $\ul{\bf H}$ is a boundary point of ${\cal
S}_I(I,I,2)$, it follows that $\ul{\bf R}$ itself must be a boundary point
of ${\cal S}_I(I,I,2)$.\ep

\noindent We are now ready to present our result for $I=J=R$.

\begin{thm}
\label{t-1}
It holds that ${\cal P}_I(I,I,2)=\overline{\cal S}_I(I,I,2)$.
\end{thm}

\noindent {\bf Proof.} First, observe that $\ul{\bf X}=({\bf I}_I,{\bf
I}_I,{\bf U})\cdot\ul{\bf Y}$ with nonsingular ${\bf U}$ has a full GSD if
and only if $\ul{\bf Y}$ has a full GSD. Indeed, a slicemix of upper
triangular slices results in upper triangular slices.

This observation, together with Stegeman and De Lathauwer \cite[Lemma
5.1]{SDL}, yields the following results for the arrays $\ul{\bf Y}$ in
Proposition~\ref{p-1}. If $\ul{\bf Y}$ satisfies (a1) or (a2), then
$\ul{\bf Y}\in{\cal P}_I(I,I,2)$. If $\ul{\bf Y}$ satisfies (a3), then
$\ul{\bf Y}\notin{\cal P}_I(I,I,2)$.

Lemma~\ref{lem-QZ} shows that an array $\ul{\bf Y}$ satisfying (b) lies in
${\cal P}_I(I,I,2)$. Since (a)-(b) defines a partition of $\R^{I\times
I\times 2}$, we have shown that $\ul{\bf Y}\in{\cal P}_I(I,I,2)$ if and
only if $\ul{\bf Y}\in\overline{\cal S}_I(I,I,2)$. This completes the
proof.\ep

\section{Extension to general $I,J,R$}
\setcounter{equation}{0}
A GSD exists only for $R\le\min(I,J)$. However, nonexistence of an optimal
solution to problem (\ref{prob-CP}) for $I\times J\times 2$ arrays does not
seem to occur for $R>I$ or $R>J$; see Stegeman \cite{Ste2}. In
Theorem~\ref{t-2} below we show that ${\cal P}_R(I,J,2)=\overline{\cal
S}_R(I,J,2)$ for $R\le\min(I,J)$. This extends Theorem~\ref{t-1}. In the
proof of Theorem~\ref{t-2}, we make use of Theorem~\ref{t-1} and the
following lemma, which concerns an orthogonal equivalence between interior
and boundary points of ${\cal S}_R(I,J,2)$ and those of ${\cal
S}_R(R,R,2)$.

\begin{lem}
\label{lem-orthequiv}
Let $\ul{\bf Y}\in\R^{I\times J\times 2}$ with $R\le\min(I,J)$. Then $\ul{\bf
Y}\in\overline{\cal S}_R(I,J,2)$ if and only if there exist
${\bf S}$ $(I\times R)$ and ${\bf T}$ $(J\times R)$ with ${\bf S}^T{\bf
S}={\bf T}^T{\bf T}={\bf I}_R$ such that $\ul{\bf Y}=({\bf S},{\bf T},{\bf
I}_2)\cdot\ul{\bf X}$ with $\ul{\bf X}\in\overline{\cal S}_R(R,R,2)$.
Moreover, $\ul{\bf Y}\in{\cal S}_R(I,J,2)$ if and only if $\ul{\bf
X}\in{\cal S}_R(R,R,2)$. \end{lem}

\noindent {\bf Proof.} See \cite[Theorem 5.2]{DSL}.\ep

\begin{thm}
\label{t-2}
Let $R\le\min(I,J)$.
It holds that ${\cal P}_R(I,J,2)=\overline{\cal S}_R(I,J,2)$.
\end{thm}

\noindent {\bf Proof.} Let $\ul{\bf Y}\in\overline{\cal S}_R(I,J,2)$. By
Lemma~\ref{lem-orthequiv} and Theorem~\ref{t-1} we have $\ul{\bf Y}=({\bf
S},{\bf T},{\bf I}_2)\cdot\ul{\bf X}$ with $\ul{\bf X}\in\overline{\cal
S}_R(R,R,2)={\cal P}_R(R,R,2)$. This implies
\be
\label{eq-need1}
{\bf Y}_k={\bf S}\,{\bf X}_k\,{\bf T}^T=({\bf S}{\bf Q}_a)\,{\bf
R}_k\,({\bf T}{\bf Q}_b)^T\,,\quad\quad k=1,2\,. \ee

\noindent Since the matrices ${\bf S}{\bf Q}_a$ and ${\bf T}{\bf Q}_b$ are
column-wise orthonormal and ${\bf R}_k$ is $R\times R$ upper triangular,
(\ref{eq-need1}) implies that $\ul{\bf Y}\in{\cal P}_R(I,J,2)$.

Next, let $\ul{\bf Y}\in{\cal P}_R(I,J,2)$. Then ${\bf Y}_k={\bf Q}_a\,{\bf
R}_k\,{\bf Q}_b^T$ for $k=1,2$, which is equivalent to $\ul{\bf Y}=({\bf
Q}_a,{\bf Q}_b,{\bf I}_2)\cdot\ul{\bf R}$, where $\ul{\bf R}\in\R^{R\times
R\times 2}$ has two upper triangular slices. Hence, $\ul{\bf R}\in{\cal
P}_R(R,R,2)=\overline{\cal S}_R(R,R,2)$ by Theorem~\ref{t-1}. An application
of Lemma~\ref{lem-orthequiv} yields $\ul{\bf Y}\in\overline{\cal
S}_R(I,J,2)$. This completes the proof.\ep

\section{Conclusion}
\setcounter{equation}{0}
We have shown that the set of $I\times J\times 2$ arrays with a full GSD of
size $R$ equals the closure of the set of $I\times J\times 2$ arrays with
at most rank $R$. Also, we have provided a complete classification of
interior, boundary, and exterior points of the latter set. This extends the
theoretical results in \cite{SDL}, which were limited to the case of
nonsingular upper triangular matrices in the GSD.

\end{document}